\documentclass[10pt]{article}

\pdfoutput=1

\usepackage[bookmarks=true]{hyperref}  

\usepackage{ifpdf}
\ifpdf
   \usepackage[pdftex]{graphicx}
   \usepackage{epstopdf}
\else
   \usepackage{graphicx}
\fi

\usepackage{color}
\usepackage{amssymb,amsmath}

\usepackage{calc}

\usepackage{tikz}

\newcommand{\xx}{x}

\newcommand{\bu}{{\bf u}}
\newcommand{\fx}{{\bf f}}
\newcommand{\ux}{u}

\newcommand{\junk}[1]{{}}

\begin{document}
\title{A Note on the Convergence of the Godunov Method for Impact Problems}

\author{J. W. Banks \\ 
\small banks20@llnl.gov\\
\small 7000 East Ave., L-422\\
\small Livermore, CA 94551}

\maketitle

\begin{abstract}
This paper identifies a new pathology that can be found for numerical simulations of nonlinear conservation law systems. Many of the difficulties already identified in the literature (rarefaction shocks, carbuncle phenomena, slowly moving shocks, wall heating, etc) can be traced to insufficient numerical dissipation, and the current case is no different. However, the details of the case we study here are somewhat unique in that the solution which is found by the numerics is very weak and can fail to have a derivative anywhere in the post-shock region.
\end{abstract}


\section{Introduction} \label{sec:introduction}
Theoretical convergence characteristics of numerical methods for systems of nonlinear PDEs have been difficult to ascertain even in 1D. With the notable exception of the random choice method of Glimm~\cite{glimm65} and its extensions, rigorous error bounds have remained largely elusive. This is not a situation where there is simply a hole in the theory but convergence troubles are never found in practice. There are in fact a number of well-known examples where numerical methods are known to behave poorly. Examples in a single space dimension include rarefaction shocks at sonic points~\cite{harten76}, the so-called wall heating phenomenon~\cite{noh87,rider00}, and sub-linear convergence for linear waves~\cite{hedstrom68,banks08a}. Often these difficulties are associated with discontinuous solutions or a lack of sufficient dissipation in the method because the nonlinear artificial dissipation inherent to the schemes goes to zero at certain points in the flow. We investigate here a new pathology where the dissipation is insufficient over a large portion of the domain. The result is convergence to a very weak solution which is nowhere differentiable. The exact nature of this weak solution is seen to depend heavily on the choice of time step. Admittedly this poor behavior could be eliminated with a simple linear artificial viscosity, but the intent here is simply to indicate the kind of behavior that can be found.

\section{Governing equations and model problem} \label{sec:model}

Consider the one-dimensional Euler equations with ideal equation of state
\begin{equation}
  {\partial\over\partial t}\bu+{\partial\over\partial\xx}\fx(\bu)=0,
  \label{eq:governing}
\end{equation}
where $\bu=[\rho, \rho\ux, \rho E,]^{T}$ and $\fx(\bu)=[\rho\ux, \rho\ux^2+p, \ux(\rho E+p)]^{T}$. Here $\rho$ is the density, $u$ the velocity, $E$ the total energy per unit mass, and $p$ the pressure. The total energy for the fluid is given by $E=e+\frac{1}{2}\ux^2$ where the equation of state is given by $e=\frac{p}{\rho(\gamma-1)}$ with $\gamma$ the ratio of specific heats. 

We investigate an impact problem with $\gamma=1.4$ on the domain $x\in[-.5,.5]$. The initial conditions in primitive variables are
\[
  \left[\rho,u,p\right]= \left\{
    \begin{array}{lc}
      \left[1.0,2.0,\frac{1}{\gamma}\right] & \qquad \hbox{for $x<0$} \smallskip \\
      \left[1.0,-2.0,\frac{1}{\gamma}\right] & \qquad \hbox{for $x\ge0$}.
    \end{array}
  \right.
\]
Inflow conditions are applied at domain boundaries, and we integrate to time $t_f=0.5$. 

\section{Numerical results} \label{sec:results}
We approximate the solution using the first-order Godunov method~\cite{godunov59} with Roe's approximate Riemann solver~\cite{roe81,toro99}. Note that the results do not change in any significant way if one instead uses an exact Riemann solver, nor if one moves to a second-order or high-resolution scheme. Discretization is performed on the computational mesh $x_i=-0.5+(i-1)\Delta x$ for $i=1\ldots m$ and $\Delta x = 1/(m-1)$.  Initial conditions are applied with exact states to the left and right of the origin. For the cases considered here, $x_i\ne0$ for any $i$ and the initial condition is applied as an exact conservative average of the left and right conservative states.

\begin{figure}
\begin{center}
  \includegraphics[width=2.3in]{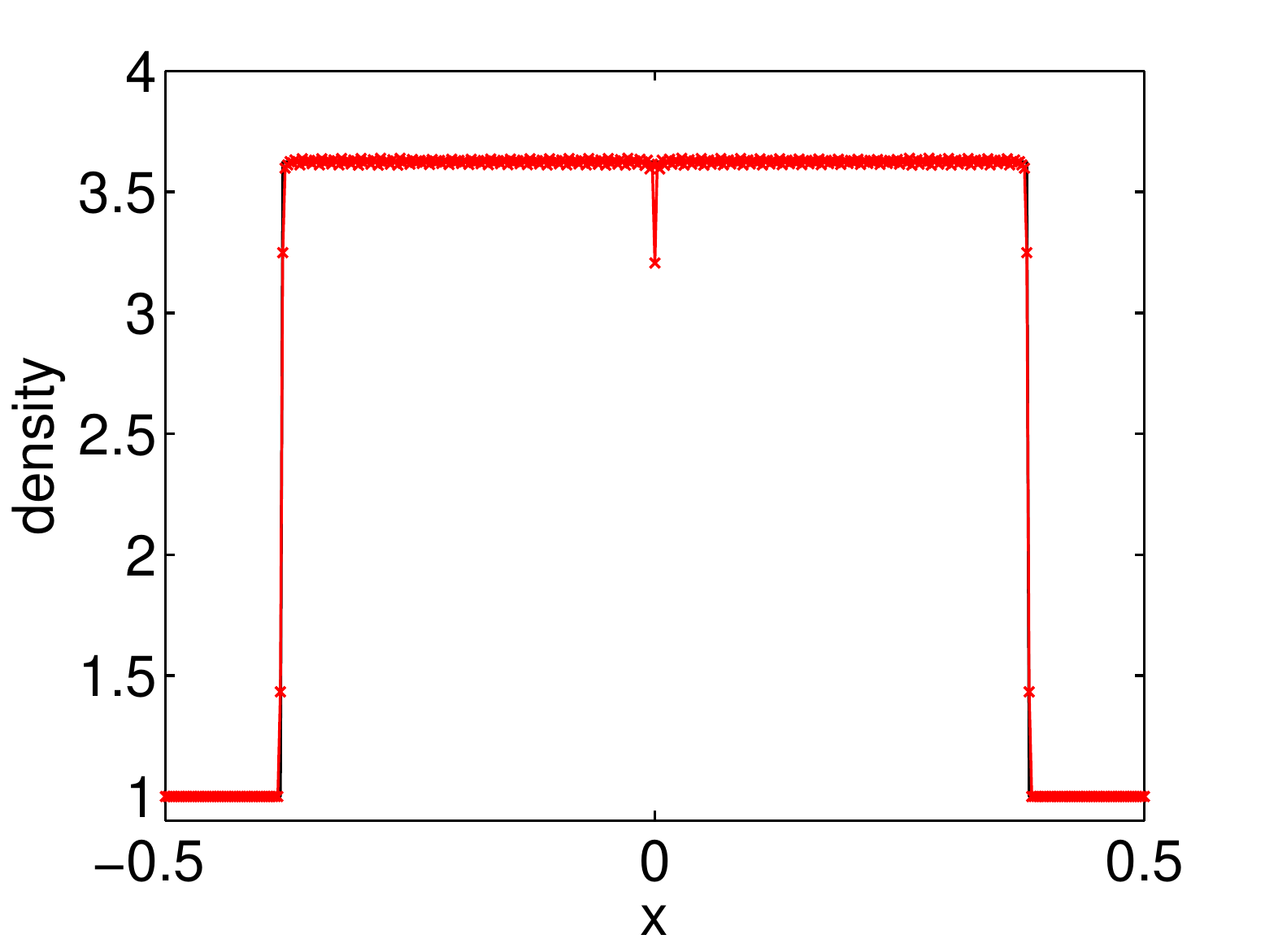} \hspace{0.1in}
  \includegraphics[width=2.3in]{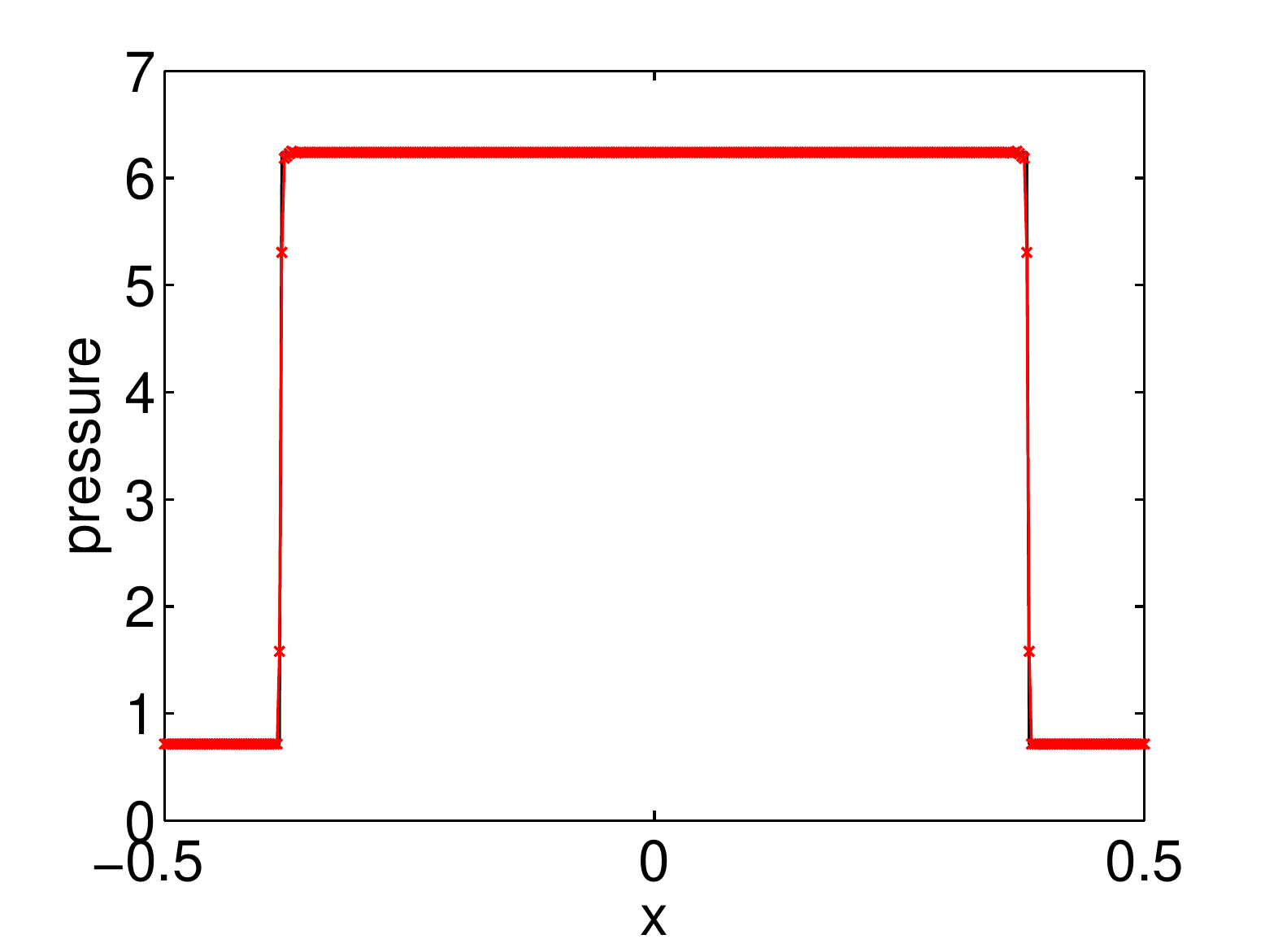} \vspace{0.1in}\\
  \includegraphics[width=2.3in]{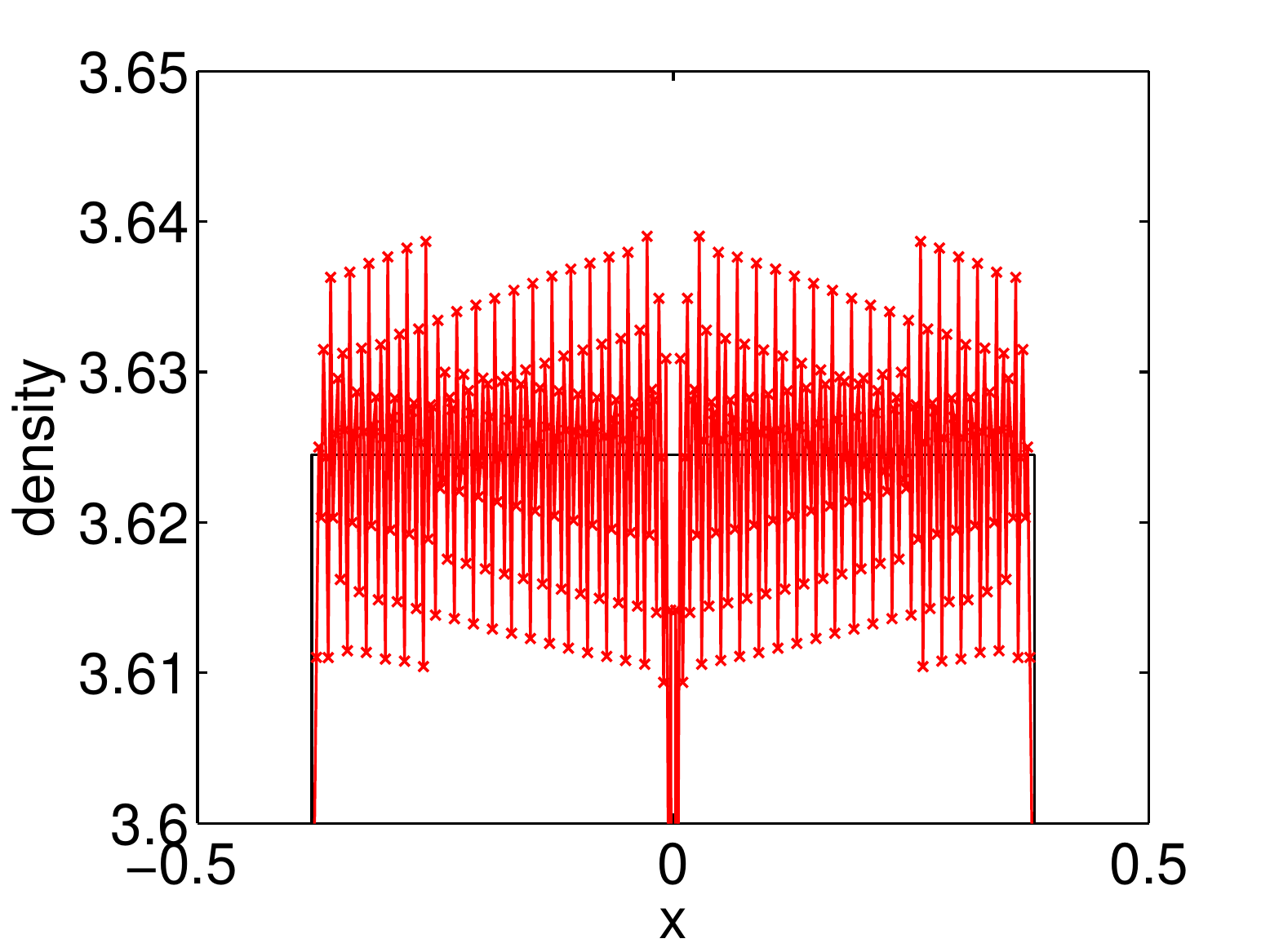} \\ 
  \caption{Density (top left), pressure (top right), and a zoom of the density (bottom). In all plots the black line is the entropy satisfying solution and the red 'x' marks and corresponding line are the first-order Godunov approximation.}
  \label{fig:impactN400}
\end{center}
\end{figure}

Approximate solutions to the impact problem are shown in Figure~\ref{fig:impactN400}.
Shown are the density, pressure, and a restricted view of the density for a CFL number of $0.9$ and $m=401$. The eventual limiting behavior can already seen in the zoom of the density where the solution oscillates around the entropy satisfying solution. Because the artificial viscosity in the Godunov method is dependent on the velocity and because the exact solution has no post-shock velocity for this problem, the magnitude of the oscillations does not decrease if the post-shock velocity converges to zero fast enough in some sense. This appears to be the case and the frequency of oscillation increases in an unbounded manner as the grid resolution increases. Such an approximation will be correct in some average sense, but will not converge in an $L_2$ or even $L_1$ sense. The solution to which the numerical approximation is converging appears to be a measure valued solution, whose value at a point can be drawn from a statistical distribution. Intuitively one can think that as $\Delta x \to 0$, the approximate solution in the post-shock region is converging to a solution that lives in a band containing the exact solution. The size of that band is dependent on the details of the discretization, and most prominently the time step. Furthermore the width can be zero for certain circumstances, as shown below.

\begin{figure}
\begin{center}
  \includegraphics[width=2.3in]{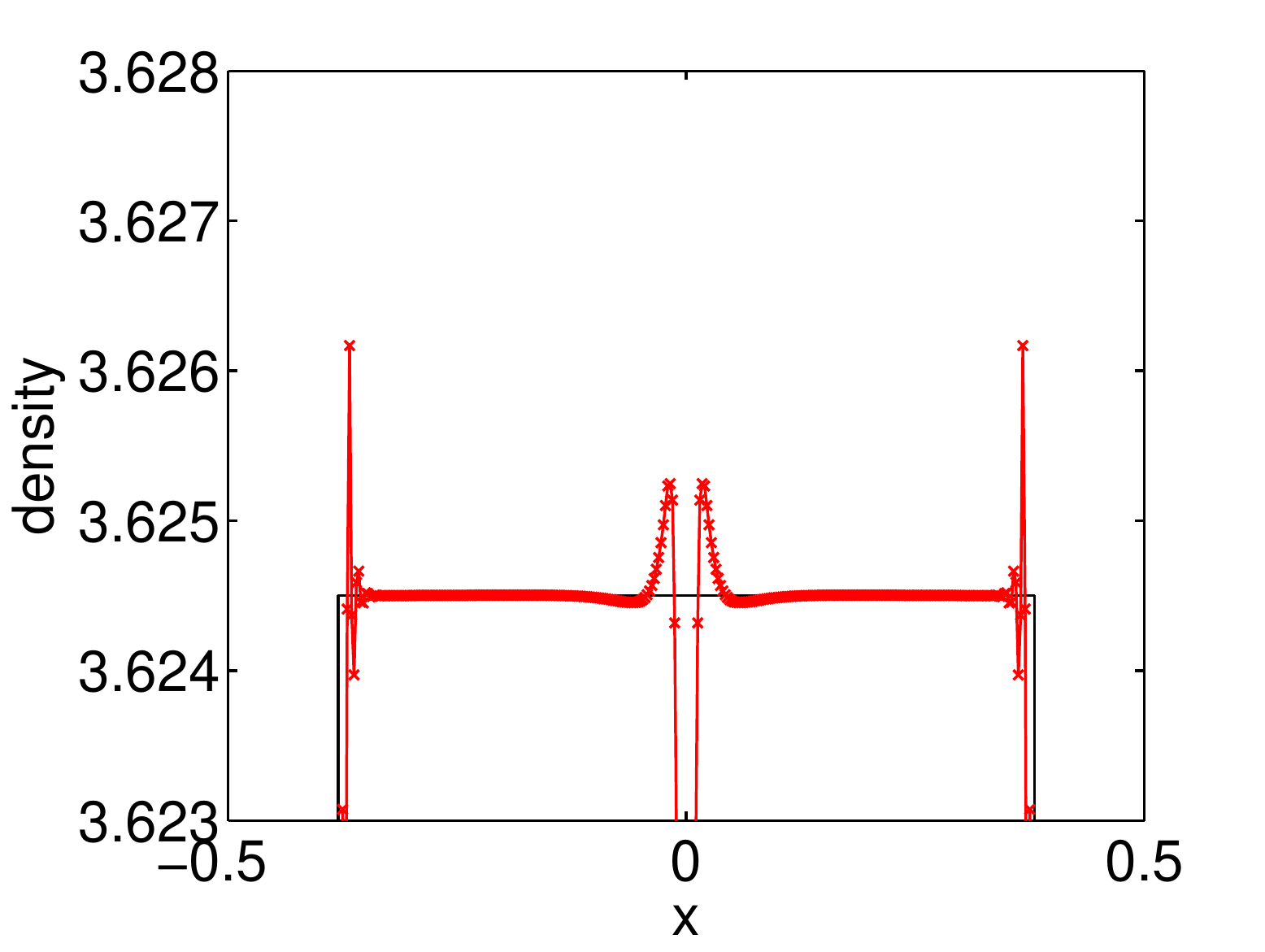} \hspace{0.1in}
  \includegraphics[width=2.3in]{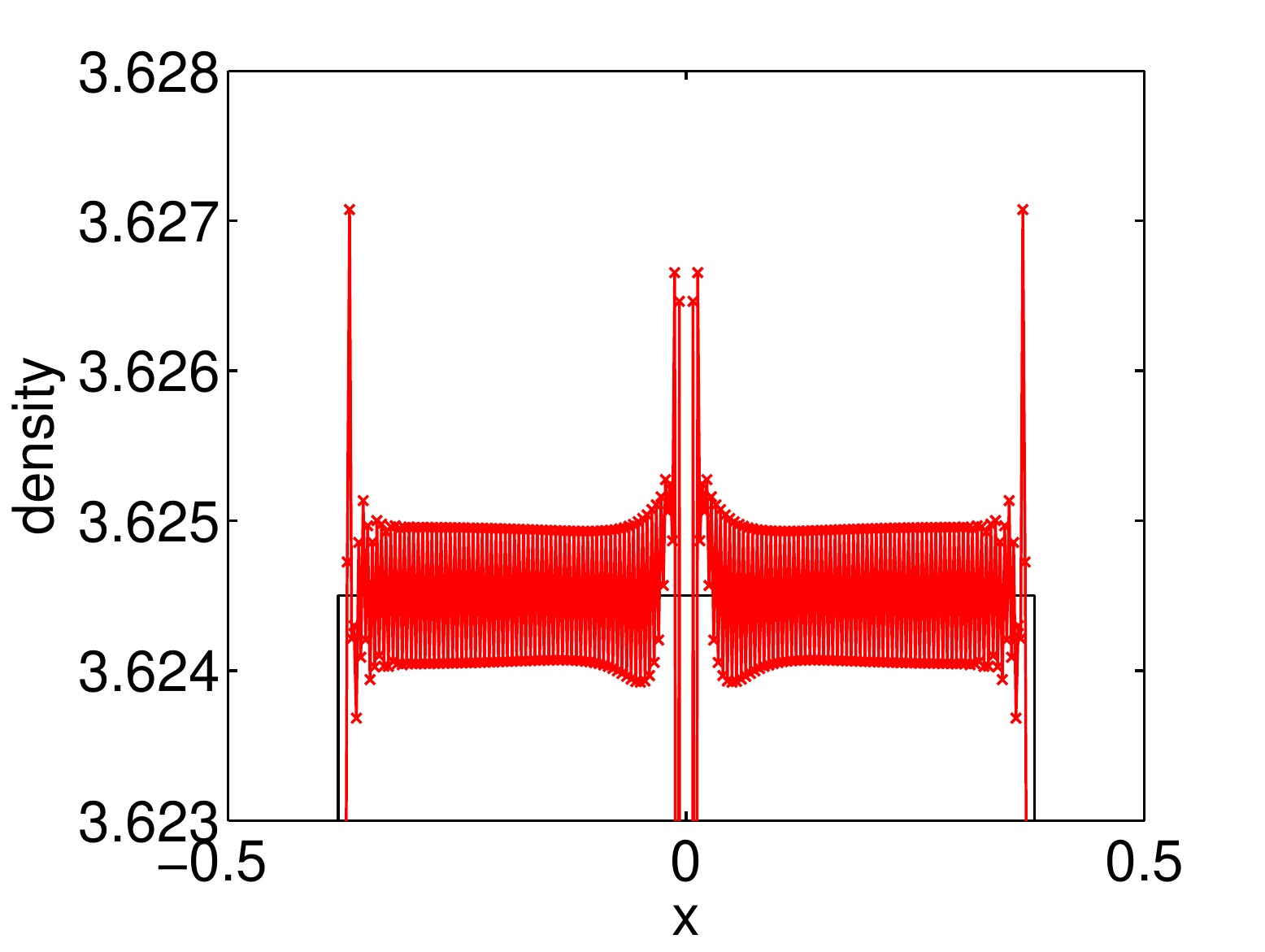}
  \caption{Zoom of the density where the shock travels through one computational cell in exactly $5$ time steps (left) and exactly $5.5$ steps (right).}
  \label{fig:impactN400_integral}
\end{center}
\end{figure}

\begin{figure}[hbt]
\begin{center}
\begin{tikzpicture}[scale=1]
    \draw(0,0) node[anchor=west] {\includegraphics[height=4cm]{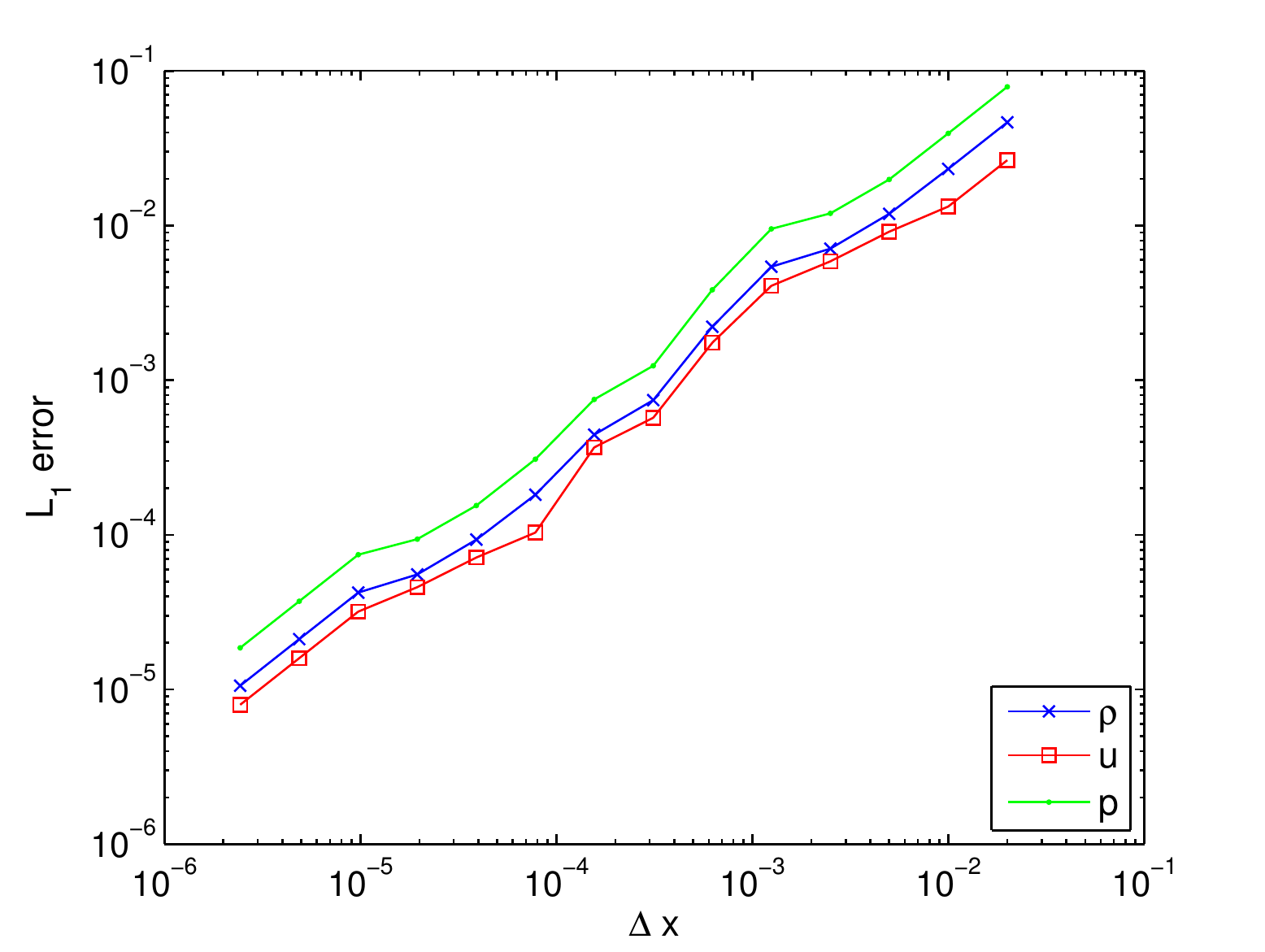}};
    \draw(5.5,0) node[anchor=west]{
{\footnotesize
\begin{tabular}{||c||c|c|c|c|c|c||} \hline\hline
  $m$ & $e_\rho(m)$ & $\kappa$ & $e_u(m)$ &$\kappa$ & $e_p(m)$ & $\kappa$ \\ \hline\hline
51 & $4.646$e$-02$ & -- & $2.649$e$-02$ & -- & $7.896$e$-02$ & -- \\ \hline 
101 & $2.325$e$-02$ & $1.00$ & $1.327$e$-02$ & $1.00$ & $3.945$e$-02$ & $1.00$ \\ \hline 
201 & $1.191$e$-02$ & $0.97$ & $9.133$e$-03$ & $0.54$ & $1.984$e$-02$ & $0.99$ \\ \hline 
401 & $7.097$e$-03$ & $0.75$ & $5.873$e$-03$ & $0.64$ & $1.201$e$-02$ & $0.72$ \\ \hline 
801 & $5.412$e$-03$ & $0.39$ & $4.078$e$-03$ & $0.53$ & $9.524$e$-03$ & $0.33$ \\ \hline 
1601 & $2.217$e$-03$ & $1.29$ & $1.746$e$-03$ & $1.22$ & $3.839$e$-03$ & $1.31$ \\ \hline 
3201 & $7.444$e$-04$ & $1.57$ & $5.706$e$-04$ & $1.61$ & $1.240$e$-03$ & $1.63$ \\ \hline 
6401 & $4.436$e$-04$ & $0.75$ & $3.671$e$-04$ & $0.64$ & $7.506$e$-04$ & $0.72$ \\ \hline 
12801 & $1.818$e$-04$ & $1.29$ & $1.036$e$-04$ & $1.82$ & $3.081$e$-04$ & $1.28$ \\ \hline 
25601 & $9.305$e$-05$ & $0.97$ & $7.133$e$-05$ & $0.54$ & $1.549$e$-04$ & $0.99$ \\ \hline 
51201 & $5.545$e$-05$ & $0.75$ & $4.589$e$-05$ & $0.64$ & $9.383$e$-05$ & $0.72$ \\ \hline
102401 & $4.228$e$-05$ & $0.29$ & $3.186$e$-05$ & $0.53$ & $7.441$e$-05$ & $0.33$ \\ \hline
204801 & $2.114$e$-05$ & $1.00$ & $1.593$e$-05$ & $1.00$ & $3.720$e$-05$ & $1.00$ \\ \hline
409601 & $1.057$e$-05$ & $1.00$ & $7.965$e$-06$ & $1.00$ & $1.860$e$-05$ & $1.00$ \\ \hline\hline
 \end{tabular}}};
 \end{tikzpicture}
\caption{$L_1$ errors and convergence rates for density, velocity, and pressure for time step chosen so that the shock travels through each cell in exactly $5$ time steps. This corresponds to CFL $\approx .787$.}
\label{table:conv_5steps}
\end{center}
\end{figure}

We perform a convergence study using the discrete $L_1$ norm to judge convergence.
\begin{figure}[hbt]
\begin{center}
\begin{tikzpicture}[scale=1]
    \draw(0,0) node[anchor=west] {\includegraphics[height=4cm]{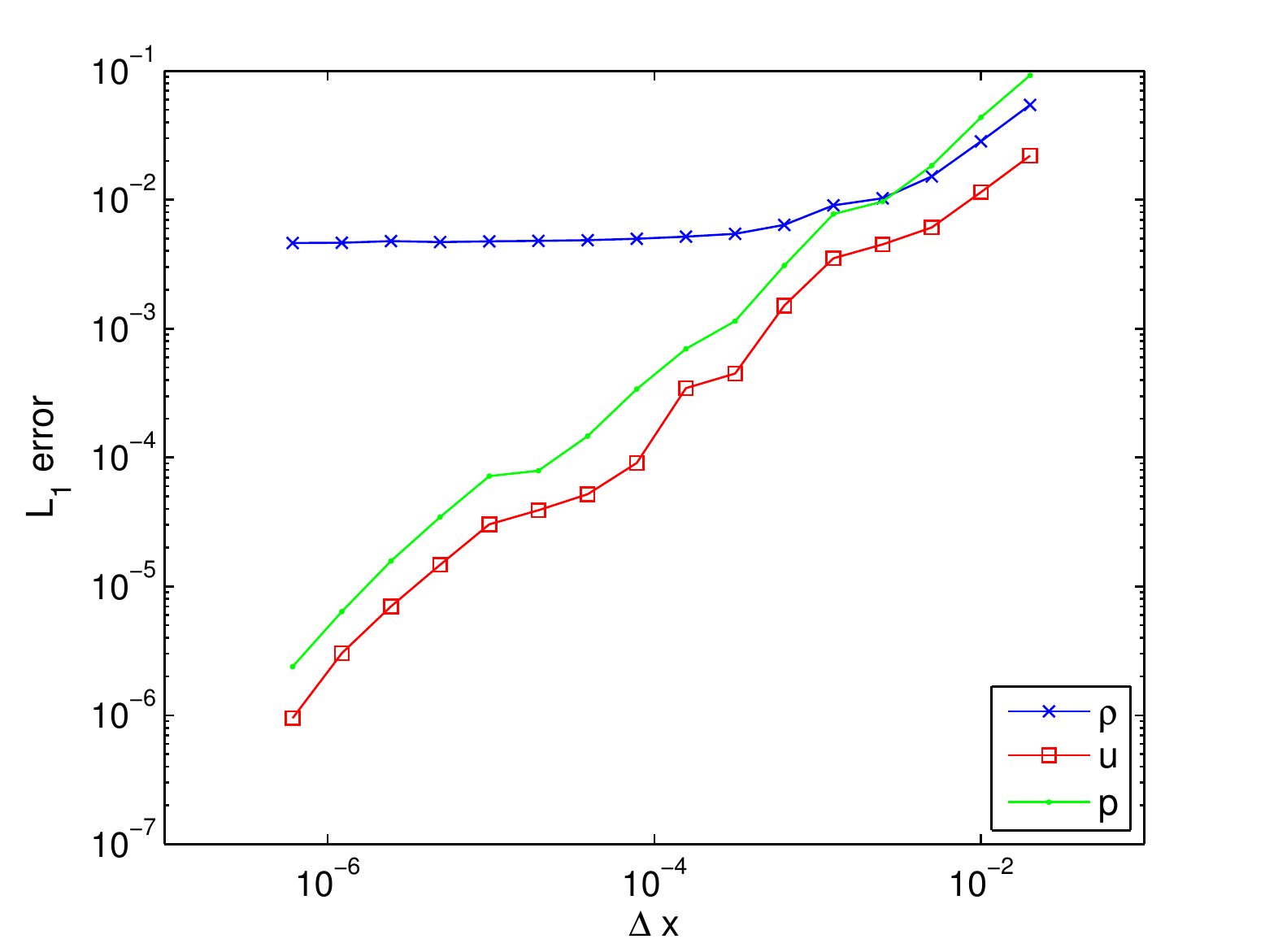}};
    \draw(5.5,0) node[anchor=west]{
{\footnotesize
\begin{tabular}{||c||c|c|c|c|c|c||} \hline\hline
  $m$ & $e_\rho(m)$ & $\kappa$ & $e_u(m)$ &$\kappa$ & $e_p(m)$ & $\kappa$ \\ \hline\hline
51 & $5.444$e$-02$ & -- & $2.199$e$-02$ & -- & $9.267$e$-02$ & -- \\ \hline
101 & $2.838$e$-02$ & $0.94$ & $1.146$e$-02$ & $0.94$ & $4.357$e$-02$ & $1.09$ \\ \hline
201 & $1.517$e$-02$ & $0.90$ & $6.098$e$-03$ & $0.91$ & $1.841$e$-02$ & $1.24$ \\ \hline
401 & $1.028$e$-02$ & $0.56$ & $4.502$e$-03$ & $0.44$ & $9.670$e$-03$ & $0.93$ \\ \hline
801 & $9.064$e$-03$ & $0.18$ & $3.516$e$-03$ & $0.36$ & $7.760$e$-03$ & $0.32$ \\ \hline
1601 & $6.375$e$-03$ & $0.51$ & $1.509$e$-03$ & $1.22$ & $3.084$e$-03$ & $1.33$ \\ \hline
3201 & $5.438$e$-03$ & $0.23$ & $4.491$e$-04$ & $1.75$ & $1.146$e$-03$ & $1.43$ \\ \hline
6401 & $5.185$e$-03$ & $0.07$ & $3.460$e$-04$ & $0.38$ & $6.985$e$-04$ & $0.71$ \\ \hline
12801 & $4.983$e$-03$ & $0.06$ & $9.070$e$-05$ & $1.93$ & $3.398$e$-04$ & $1.04$ \\ \hline
25601 & $4.864$e$-03$ & $0.03$ & $5.198$e$-05$ & $0.80$ & $1.468$e$-04$ & $1.21$ \\ \hline
51201 & $4.804$e$-03$ & $0.02$ & $3.907$e$-05$ & $0.41$ & $7.907$e$-05$ & $0.89$ \\ \hline
102401 & $4.764$e$-03$ & $0.01$ & $3.029$e$-05$ & $0.37$ & $7.161$e$-05$ & $0.14$ \\ \hline
204801 & $4.707$e$-03$ & $0.02$ & $1.473$e$-05$ & $1.04$ & $3.452$e$-05$ & $1.05$ \\ \hline
409601 & $4.778$e$-03$ & $-0.02$ & $6.984$e$-06$ & $1.17$ & $1.575$e$-05$ & $1.13$ \\ \hline
819201 & $4.644$e$-03$ & $0.04$ & $3.031$e$-06$ & $1.20$ & $6.367$e$-06$ & $1.31$ \\ \hline
1638401 & $4.628$e$-03$ & $0.005$ & $9.497$e$-07$ & $1.67$ & $2.387$e$-06$ & $1.42$ \\ \hline\hline
 \end{tabular}}};
 \end{tikzpicture}
\caption{$L_1$ errors and convergence rates for density, velocity, and pressure for CFL$=0.9$.}
\label{table:conv_CFLp9}
\end{center}
\end{figure}
The results are shown in the table of Figure~\ref{table:conv_CFLp9}. We can see that the density does not converge below approximately $5\times 10^{-3}$. To understand the details, consider a single cell as it transitions from before to after the shock. Because the method is conservative and seems to be converging to some weak solution, the approximate shock location should be correct in the limit. However, each point travels through the shock in a slightly different manner which leads to the oscillations in density. To see this more clearly we can modify the time step such that the shock travels through a computational cell in an integral number of time steps. It is determined that the shock speed is $\approx .76205$. For the results in Figure~\ref{fig:impactN400} (with CFL $=0.9$ and $m=401$), the time step is found to be $\Delta t\approx 7.496\times 10^{-4}$. As a result, the shock travels through each cell in approximately $4.3765$ time steps. For the case of $m=401$, a time step of $\Delta t\approx 6.5613\times 10^{-4}$ corresponding to CFL $\approx .787$, will have the shock traveling through each computational cell in exactly $5$ time steps. Figure~\ref{fig:impactN400_integral} demonstrates that this choice does indeed remove the post-shock density oscillations.
Also in this figure we show the results for $\Delta t=5.965\times 10^{-4}$ corresponding to CFL $\approx .716$ where the shock travels through one cell in exactly $5.5$ time steps. Here we see that every other cell travels through the shock in the same way with the result that the density oscillates cell to cell between two values. The table in Figure~\ref{table:conv_5steps} shows the results of a grid convergence study with CFL $\approx.787$ where the shock travels through each cell in exactly 5 time steps. Here the density is seen to converge well even to the finest resolution.

\section{Conclusions} \label{sec:conclusions}
The pathology of convergence of numerical approximations for impact problems studied here is an interesting one. We showed that depending on the time step, the weak solution found by the Godunov method can be convergent only in some average sense. The density in the post-shock region is found to depend on the number of time steps required for a shock to cross a single computational cell. By choosing the time step so that the shock crosses a cell in an integral number of time steps, pointwise convergence can be obtained everywhere except a set with zero measure. Other values of the time step result in oscillatory behavior which causes the derivative to not exist everywhere behind the shock. 

\section*{Acknowledgements}
\noindent
This work was performed under the auspices of the
U.S. Department of Energy by Lawrence Livermore National Laboratory
under contract number DE-AC52-07NA27344.  

\bibliographystyle{elsart-num}
\bibliography{../journal-ISI,../jwb}

\end{document}